# A new Poisson-type deviation inequality for Markov jump processes with positive Wasserstein curvature

ALDÉRIC JOULIN

*Université de Toulouse, INSA, Institut de Mathématiques de Toulouse, F-31077 Toulouse, France. E-mail: ajoulin@insa-toulouse.fr*

The purpose of this paper is to extend the investigation of Poisson-type deviation inequalities started by Joulin (*Bernoulli* **13** (2007) 782–798) to the empirical mean of positively curved Markov jump processes. In particular, our main result generalizes the tail estimates given by Lezaud (*Ann. Appl. Probab.* **8** (1998) 849–867, *ESAIM Probab. Statist.* **5** (2001) 183–201). An application to birth–death processes completes this work.

## 1. Introduction

Let $(X_t)_{t\geq 0}$ be an ergodic Markov process on a Polish state space $\mathcal{X}$, with stationary distribution $\pi$. The well-known ergodic theorem asserts that for any integrable function $\phi \in L^1(\pi)$, the empirical mean $t^{-1}\int_0^t \phi(X_s)\,\mathrm{d}s$ converges in probability to the average $\pi(\phi) := \int_{\mathcal{X}} \phi\,\mathrm{d}\pi$ as $t$ goes to infinity. Although large deviations theory gives the speed of convergence at infinity, such an asymptotic bound is unsatisfactory when one wants to estimate the minimum time to run the simulation algorithm in order to achieve a prescribed level of accuracy. Actually, the problem of finding non-asymptotic estimates has been raised and addressed by several authors. Using the Lumer–Philips theorem for a general Markov process $(X_t)_{t\geq 0}$, Wu (2000) derived an exponential decay on the deviation probability

$$\mathbb{P}\left(\left|\frac{1}{t}\int_0^t \phi(X_s)\,\mathrm{d}s - \pi(\phi)\right| \geq y\right), \qquad y > 0, \tag{1.1}$$

available for any fixed time $t$. Although Wu's estimate is sharp in large time, such an upper bound is not explicit in the parameter $y$. More recently, this result has been extended in the diffusion framework by various authors who obtained qualitative upper bounds on (1.1), provided the stationary distribution $\pi$ satisfies some functional inequalities such as Poincaré, log-Sobolev or transportation-type inequalities; see, for instance, the recent articles of Cattiaux and Guillin (2008), Djellout *et al.* (2004), Gourcy and Wu (2006) or Guillin *et al.* (2009). However, the functional inequalities approach does not seem to be relevant for Markov jump processes because this theory is not yet well developed for discrete gradients. To the author's knowledge, the problem of determining non-asymptotic upper bounds on the deviation probability (1.1) in this context has been investigated by few authors. For instance, under a spectral gap assumption and using Kato's perturbation theory for linear operators, Lezaud (1998, 2001) established Poisson-type deviation bounds, that is, upper bounds of the order $e^{-ty \log(y)}$ for large $y$, provided the function $\phi$ and the generator of the process are bounded. On the other hand, in the case of birth–death processes admitting a so-called Lipschitz spectral gap, Liu and Ma (2009) recently extended such tail estimates to Lipschitz functions $\phi$ by using martingale techniques and convex concentration inequalities.

The purpose of this paper is to present a new Poisson-type upper bound for the deviation probability (1.1) for a general Markov jump process $(X_t)_{t\geq 0}$. Our approach relies on the notion of Wasserstein curvature recently investigated by Joulin (2007), where several tail estimates were obtained for the random variable $\phi(X_t)$. Hence we extend in this article our previous work to the path-dependent integral $t^{-1} \int_0^t \phi(X_s)\,ds$. In essence, the Wasserstein curvature characterizes a contraction property of the associated semigroup on the space of probability measures on $\mathcal{X}$, endowed with a suitable Wasserstein distance. Since the positively curved case is closely related to the speed of ergodicity of the process, we expect to obtain under this assumption a convenient upper bound on (1.1) in large time.

The paper is organized as follows: in Section 2, we recall the definition of the Wasserstein curvature of a Markov jump process $(X_t)_{t\geq 0}$. Next, we state the main contribution of the paper, Theorem 2.6, in which a Poisson-type deviation bound is established in the positively curved case for the empirical mean $t^{-1} \int_0^t \phi(X_s)\,ds$, where $\phi$ is only Lipschitz. Hence we extend the tail estimates given in the bounded case by Lezaud (1998, 2001). Section 3 is devoted to the proof of Theorem 2.6, which is rather technical and divided into several lemmas. The key point of the proof corresponds to Lemma 3.2, with the tensorization of a Laplace transform. Section 4 is devoted to the case of birth–death processes. More precisely, we compute the explicit expression of the Wasserstein curvature with respect to a large class of metrics on $\mathbb{N}$. In particular, by choosing a convenient metric related to the transition rates of the associated generator, we are able to apply our deviation inequality to birth–death processes with non-necessarily bounded generator such as the classical $M/M/\infty$ queueing process.

## 2. Preliminaries and main result

Throughout the paper, $\mathcal{X}$ is a Polish space endowed with a metric $d$, the space $\mathscr{B}(\mathcal{X})$ consists of bounded measurable functions on $\mathcal{X}$ equipped with the supremum norm $\|f\|_\infty = \sup_{x \in \mathcal{X}} |f(x)|$ and $\mathrm{Lip}_d(\mathcal{X})$ is the space of Lipschitz functions on $\mathcal{X}$ with a Lipschitz seminorm defined by

$$\|f\|_{\mathrm{Lip}_d} := \sup_{x \neq y} \frac{|f(x) - f(y)|}{d(x, y)} < +\infty.$$

On a filtered probability space $(\Omega, \mathscr{F}, (\mathscr{F}_t)_{t \geq 0}, \mathbb{P})$, let $\{(X_t)_{t \geq 0}, (\mathbb{P}_x)_{x \in \mathcal{X}}\}$ be an $\mathcal{X}$-valued cadlag Markov jump process with a generator given for any function $f \in \mathscr{B}(\mathcal{X})$ by

$$\mathcal{L}f(x) = \int_\mathcal{X} (f(y) - f(x)) Q(x, \mathrm{d}y), \qquad x \in \mathcal{X}.$$

Here the transition kernel $Q$ is assumed to be stable and conservative: for any $x \in \mathcal{X}$ and any Borel set $A$,

$$Q(x, \mathcal{X}) < +\infty, \qquad \lim_{t \downarrow 0} \frac{P_t(x, A) - 1_A(x)}{t} = Q(x, A) - Q(x, \mathcal{X}) 1_A(x),$$

where $P_t(x, \mathrm{d}y) := \mathbb{P}_x(X_t \in \mathrm{d}y)$ denotes the transition probability of the process. Let $(P_t)_{t \geq 0}$ be the associated Markov semigroup acting on the space $\mathscr{B}(\mathcal{X})$ as follows:

$$P_t f(x) := \mathbb{E}_x[f(X_t)] = \int_\mathcal{X} f(y) P_t(x, \mathrm{d}y), \qquad x \in \mathcal{X}.$$

Denote by $\mathscr{P}_d(\mathcal{X})$ the space of probability measures $\mu$ on $\mathcal{X}$ such that $\int_\mathcal{X} d(x, y) \mu(\mathrm{d}y) < +\infty$ for some (or equivalently for all) $x \in \mathcal{X}$. If the Markov kernel $P_t(x, \cdot) \in \mathscr{P}_d(\mathcal{X})$ for any $t > 0$ and any $x \in \mathcal{X}$, then the semigroup is well defined on the space $\mathrm{Lip}_d(\mathcal{X})$ and we introduce in this case the function

$$\bar{\sigma}_d(t) := -\sup\{\log \|P_t f\|_{\mathrm{Lip}_d} : \|f\|_{\mathrm{Lip}_d} = 1\}, \qquad t \geq 0,$$

with $\bar{\sigma}_d(0) = 0$. By the Markov property, the function $\bar{\sigma}_d$ is super-additive so that the following limit is well defined:

$$\sigma_d := \lim_{t \downarrow 0} \frac{\bar{\sigma}_d(t)}{t} = \inf_{t > 0} \frac{\bar{\sigma}_d(t)}{t}. \tag{2.1}$$

In particular, the number $\sigma_d$ is the best (maximal) constant $\alpha$ in the contraction inequality

$$\|P_t f\|_{\mathrm{Lip}_d} \leq e^{-\alpha t} \|f\|_{\mathrm{Lip}_d}, \qquad f \in \mathrm{Lip}_d(\mathcal{X}), \ t > 0. \tag{2.2}$$

Let us recall the definition of Wasserstein curvature of the Markov jump process $(X_t)_{t \geq 0}$ given by Joulin (2007), up to a slight modification.

**Definition 2.1.** *Assume $P_t(x, \cdot) \in \mathscr{P}_d(\mathcal{X})$ for any $t > 0$ and any $x \in \mathcal{X}$. The number $\sigma_d$ given by (2.1) is called the Wasserstein curvature of the process $(X_t)_{t \geq 0}$ with respect to the metric $d$.*

***Remark 2.2.*** In the remainder of this paper, we will remove the metric symbol $d$ in the definition of the Wasserstein curvature $\sigma_d$ when there is no risk of confusion. Moreover, we will assume implicitly that the Markov kernel $P_t(x,\cdot)$ belongs to the space $\mathscr{P}_d(\mathcal{X})$ for any $t>0$ and any $x\in\mathcal{X}$.

We define the Wasserstein distance $W_d(\mu,\nu)$ between two probability measures $\mu,\nu\in\mathscr{P}_d(\mathcal{X})$ as

$$W_d(\mu,\nu):=\inf_\gamma \int_{\mathcal{X}\times\mathcal{X}} d(x,y)\gamma(\mathrm{d}x,\mathrm{d}y),$$

where the infimum is taken over all $\gamma\in\mathscr{P}_d(\mathcal{X}\times\mathcal{X})$ with marginals $\mu$ and $\nu$. The Kantorovich–Rubinstein duality theorem allows us to rewrite the Wasserstein distance as

$$W_d(\mu,\nu)=\sup\left\{\left|\int_\mathcal{X} f\,\mathrm{d}\mu - \int_\mathcal{X} f\,\mathrm{d}\nu\right|:\|f\|_{\mathrm{Lip}_d}\leq 1\right\},$$

see, for instance, Chen (2004), Theorem 5.10. Hence the Wasserstein curvature $\sigma$ is also the best (maximal) constant $\alpha$ in the inequality

$$W_d(P_t(x,\cdot),P_t(y,\cdot))\leq \mathrm{e}^{-\alpha t}d(x,y),\qquad x,y\in\mathcal{X},\ t>0. \tag{2.3}$$

***Remark 2.3.*** As noted by Joulin (2007), our definition of Wasserstein curvature of Markov processes is inspired by the continuous setting of Brownian motion on Riemannian manifolds studied by Sturm and Von Renesse (2005), where it is stated that the contraction inequality (2.3) characterizes uniform lower bounds on the Ricci curvature of the manifold. However, after our paper was published, we learned that a similar notion of curvature for Markov processes relying on such an inequality had been previously introduced in the PhD thesis of Sammer (2005) under the name "Ricci–Wasserstein curvature", and later independently by Ollivier (2009, 2007b) as the "Ricci curvature" of Markov chains on metric spaces. Actually, without the link to geometry, the inequality (2.3) appeared first in the work of Dobrushin (1970) with his study on random fields, and is known in statistical mechanics as the "Dobrushin uniqueness condition". Moreover, such a contraction inequality is fundamental to estimate the spectral gap $\lambda_1$ (say) of reversible Markov processes, or equivalently to establish a Poincaré inequality for the stationary distribution, since we have $\lambda_1 \geq \sigma$. See, for instance, Chen (2004), Chapter 9, for a summary of and precise references for this topic.

Actually, the Wasserstein curvature is closely related to the ergodicity of the process, as illustrated by the following result. See, for instance, the very general result of Dobrushin (1970), Theorem 3, for a proof in the discrete-time case or Chen (2004), Theorem 5.23, in the continuous-time setting.

**Theorem 2.4.** *Assume $\sigma > 0$. Then the process $(X_t)_{t \geq 0}$ admits a unique stationary distribution $\pi \in \mathscr{P}_d(\mathcal{X})$ and is ergodic in the following sense: For any initial point $x \in \mathcal{X}$,*

$$W_d(P_t(x, \cdot), \pi) \leq e^{-\sigma t} \int_{\mathcal{X}} d(x, y) \pi(\mathrm{d}y) \underset{t \to +\infty}{\longrightarrow} 0. \tag{2.4}$$

**Remark 2.5.** When $d$ is the trivial metric on $\mathcal{X}$ defined by $d(x, y) = 1_{\{x \neq y\}}$, the Wasserstein distance is nothing but half of the total variation norm. Therefore, the convergence in Wasserstein distance generalizes the classical convergence in total variation used in the context of general Markov processes.

Under the ergodic property of the process, the celebrated ergodic theorem states that for any $\phi \in L^1(\pi)$, the empirical mean $t^{-1} \int_0^t \phi(X_s) \, \mathrm{d}s$ converges in probability as $t$ goes to infinity to the equilibrium $\pi(\phi) := \int_{\mathcal{X}} \phi \, \mathrm{d}\pi$, where $\pi$ denotes the unique stationary distribution given by Theorem 2.4. It is well known that the determination of qualitative non-asymptotic deviation inequalities is of fundamental importance for simulation algorithms. However, the theory of large deviations provides a bound for this convergence that is only asymptotic in time on the one hand, and whose behaviour in terms of the deviation level is not explicit on the other hand. Hence one may wonder if Wasserstein curvature plays a crucial role in the determination of such tail estimates relating the speed of ergodicity of the process. We give now an affirmative answer to this question by stating the main result of the paper, the proof of which is given in the next section. In the remainder of the paper, we denote the function

$$g(u) := (1 + u) \log(1 + u) - u, \qquad u > 0. \tag{2.5}$$

**Theorem 2.6.** *Assume $\sigma > 0$ and that there exist two positive constants $b$ and $V$ such that*

$$\sup_{t > 0} d(X_{t-}, X_t) \leq b \quad \text{and} \quad \left\| \int_{\mathcal{X}} d(\cdot, y)^2 Q(\cdot, \mathrm{d}y) \right\|_{\infty} \leq V^2. \tag{2.6}$$

*Letting $\phi \in \mathrm{Lip}_d(\mathcal{X})$, for any initial state $x \in \mathcal{X}$, any $t > 0$ and any $y > 0$ we have the Poisson-type deviation inequality:*

$$\mathbb{P}_x \left( \left| \frac{1}{t} \int_0^t \phi(X_s) \, \mathrm{d}s - \pi(\phi) \right| \geq y + M_t^x \right) \leq 2 e^{-(V^2 t / b^2) g((by\sigma)/(V^2(1 - e^{-\sigma t}) \|\phi\|_{\mathrm{Lip}_d}))}, \tag{2.7}$$

*where $\pi$ denotes the unique stationary distribution given in Theorem 2.4 and*

$$M_t^x := \frac{(1 - e^{-\sigma t}) \|\phi\|_{\mathrm{Lip}_d}}{\sigma t} \int_{\mathcal{X}} d(x, z) \pi(\mathrm{d}z) \underset{t \to +\infty}{\longrightarrow} 0.$$

Let us give some comments on this result.

**Remark 2.7.** According to a classical large deviation result, the estimate (2.7) is optimal in large time since the order of magnitude is $e^{-ct}$, and is also sharp in small time.

Moreover, the function $u \mapsto g(u)$ is equivalent to $u^2/2$ as $u$ is close to 0 and to $u \log(u)$ as $u$ tends to infinity. Hence, for sufficiently large $t$ the inequality (2.7) exhibits a Gaussian regime for small values of the deviation level $y$, in accordance with the central limit theorem for Markov processes and a Poisson regime for its large values.

*Remark 2.8.* Assume that the process is reversible. As noted in Remark 2.3, the positivity of the Wasserstein curvature ensures the existence of a spectral gap $\lambda_1$ of the underlying generator, that is, $\lambda_1 \geq \sigma > 0$. Therefore, using the Poincaré inequality, the asymptotic variance of the empirical mean is bounded by $V^2 \|\phi\|^2_{\text{Lip}_d}/\lambda_1^2$ and one deduces that the right-hand side of (2.7) is sharp in $\sigma$ in the Gaussian regime since it behaves as $e^{-t\sigma^2 y^2/(2V^2 \|\phi\|^2_{\text{Lip}_d})}$ for large time.

*Remark 2.9.* Up to constant factors, we extend the Chernoff inequalities established by Lezaud (1998, 2001), because boundedness assumptions are required neither on the function $\phi$ nor on the generator. Note, however, that if the metric $d$ is such that $\inf_{x \neq y} d(x,y) > 0$, then the finiteness of $V$ implies that the generator is bounded. In particular, when $d$ is the trivial metric, we recover Lezaud's results since we have in this case $\text{Lip}_d(\mathcal{X}) = \mathcal{B}(\mathcal{X})$ and $V^2 = \|Q(\cdot, \mathcal{X})\|_\infty$. Nevertheless, the price to pay in Theorem 2.6 is to assume $\sigma > 0$, which is a stronger assumption in the reversible case than the existence of a spectral gap required by Lezaud.

*Remark 2.10.* Consider for instance the Langevin-type diffusion process solution of the following stochastic differential equation $dX_t = \sqrt{2} \, dB_t - \nabla U(X_t) \, dt$, where $(B_t)_{t \geq 0}$ is a standard Brownian motion on the Euclidean space $(\mathbb{R}^n, d)$ and $U$ is a regular potential such that $\int e^{-U(x)} \, dx = 1$. Denote by $\pi(dx) = e^{-U(x)} \, dx$ the stationary distribution of the process $(X_t)_{t \geq 0}$. Since the Wasserstein curvature can be defined in the diffusion framework, a step-by-step adaptation of the proof of Theorem 2.6 below – especially the proof of Lemma 3.1 – entails for any Lipschitz function $\phi$ on $(\mathbb{R}^n, d)$ a Gaussian deviation inequality of the form

$$\mathbb{P}_x \left( \left| \frac{1}{t} \int_0^t \phi(X_s) \, ds - \pi(\phi) \right| \geq y + M_t^x \right) \leq 2 e^{-ty^2 \sigma^2/(2(1-e^{-\sigma t})^2 \|\phi\|^2_{\text{Lip}_d})},$$

provided the Wasserstein curvature of the process $(X_t)_{t \geq 0}$ is positive. A sufficient condition ensuring this positivity is given by the Bakry–Émery curvature criterion, see Bakry and Émery (1985), under which the authors established a logarithmic Sobolev inequality for the stationary distribution $\pi$. On the other hand, it is classical that such a functional inequality entails a similar Gaussian decay to that given above; see, for instance, Wu (2000) or the recent article of Guillin *et al.* (2009). Hence we give under comparable assumptions another proof of this Gaussian tail estimate.

*Remark 2.11.* As illustrated for birth–death processes in Section 4, it is sufficient to carry the analysis in the one-dimensional case since the Wasserstein curvature tensorizes on product spaces equipped with the $\ell^1$-metric. Indeed, for each $i = 1, \ldots, N$, consider the Markov process $(X^i_t)_{t \geq 0}$ with kernel transition $Q^i$, stationary distribution $\pi^i$ and Wasserstein curvature $\sigma^i$, all valued in the same Polish space $(\mathcal{Y}, \rho)$ to simplify. We construct the multidimensional Markov process $(X_t)_{t \geq 0}$ valued in $(\mathcal{X}, d)$, where $\mathcal{X} := \mathcal{Y}^N$ and $d$ is the $\ell^1$-metric defined with respect to $\rho$, as follows: choose first a coordinate uniformly at random and then let the univariate dynamics run according to this direction. Then the stationary distribution $\pi$ is given by $\pi = \bigotimes_{i=1}^N \pi^i$. Now let $\mu$ and $\nu$ be two product probability measures on $\mathcal{X}$. Then the classical tensorization property of the Wasserstein distance is given by $W_d(\mu, \nu) = \sum_{i=1}^N W_\rho(\mu^i, \nu^i)$, see for instance Sammer (2005), Lemma 2.2.6, for a proof. Hence, the Wasserstein curvature $\sigma$ with respect to the metric $d$ of the Markov process $(X_t)_{t \geq 0}$ is computed as $\sigma = \min_{i=1,\ldots,N} \sigma^i/N$. Moreover, if we denote by $b_i$ and $V_i$ the numbers in (2.6) related to the coordinate process $(X^i_t)_{t \geq 0}$, then Theorem 2.6 applies for the multidimensional Markov process $(X_t)_{t \geq 0}$ with $\sigma$ and $\pi$ as above and with $b := \max_{i=1,\ldots,N} b_i$ and $V^2 := \sum_{i=1}^N V_i^2/N$.

To illustrate our argument, consider the symmetric continuous-time random walk $(X_t)_{t \geq 0}$ on the discrete cube $\{0, 1\}^N$, equipped with the Hamming metric $d(x, y) = \sum_{i=1}^N 1_{\{x_i \neq y_i\}}$. The associated semigroup kernel is given by

$$P_t(x, y) = \frac{1}{2^N} \prod_{i=1}^N (1 + (-1)^{|x_i - y_i|} e^{-t/N}), \qquad x, y \in \{0, 1\}^N,$$

and the stationary distribution is the uniform probability measure on $\{0, 1\}^N$, say $\pi^{\otimes N}$. Since in the one-dimensional case a simple calculation shows that the Wasserstein curvature with respect to the trivial metric equals 1, the Wasserstein curvature on the product space with respect to the Hamming metric is $\sigma = 1/N$. Moreover, we have $b = 1$ and $V^2 = 1/2$ so that by Theorem 2.6 the following deviation inequality holds for any Lipschitz function $\phi$ with respect to the Hamming metric on $\{0, 1\}^N$:

$$\mathbb{P}_x \left( \left| \frac{1}{t} \int_0^t \phi(X_s) \, ds - \pi^{\otimes N}(\phi) \right| \geq y + M_t^x \right) \leq 2 e^{-(t/2)g(2y/(N(1-e^{-t/N})\|\phi\|_{\text{Lip}_d}))}.$$

## 3. Proof of Theorem 2.6

This section is devoted to the proof of Theorem 2.6, which is rather technical and divided into several lemmas. First, we give a convenient upper bound in large time on a univariate Laplace transform, see Lemma 3.1 below. Using the method of tensorization, the extension to the multidimensional case is considered in Lemma 3.2. Finally, with the help of the previous lemmas and by a suitable approximation of the empirical mean, we finish the proof of Theorem 2.6.

Let us establish first an upper bound on the Laplace transform of a Lipschitz function of the process $(X_t)_{t \geq 0}$. The proof, which is a straightforward adaptation of Joulin (2007), Theorem 3.1, is given for completeness.

**Lemma 3.1.** *Under the assumptions of Theorem 2.6, for any $f \in \mathrm{Lip}_d(\mathcal{X})$, any $x \in \mathcal{X}$, any $t > 0$ and any $\tau > 0$, we have*

$$\mathbb{E}_x[\mathrm{e}^{\tau(f(X_t) - \mathbb{E}_x[f(X_t)])}] \leq \exp\{h(\tau, t, b\|f\|_{\mathrm{Lip}_d})\}, \tag{3.1}$$

*where $h$ is the function defined on $(\mathbb{R}_+)^3$ by*

$$h(\tau, t, z) := \frac{V^2(1 - \mathrm{e}^{-2\sigma t})}{2b^2 \sigma}(\mathrm{e}^{\tau z} - \tau z - 1). \tag{3.2}$$

**Proof.** Assume first that the Lipschitz function $f$ is bounded. Then the process $(Z_s^f)_{0 \leq s \leq t}$ given by $Z_s^f := P_{t-s}f(X_s) - P_t f(X_0)$ is a real-valued $\mathbb{P}_x$-martingale with respect to the filtration $(\mathscr{F}_s)_{0 \leq s \leq t}$. Using (2.2) and (2.6), we have

$$\sup_{0 < s \leq t} |Z_s^f - Z_{s-}^f| = \sup_{0 < s \leq t} |P_{t-s}f(X_s) - P_{t-s}f(X_{s-})|,$$

$$\leq b\|f\|_{\mathrm{Lip}_d}$$

and also

$$\langle Z^f, Z^f \rangle_s = \int_0^s \int_{\mathcal{X}} (P_{t-\tau}f(y) - P_{t-\tau}f(X_{\tau-}))^2 Q(X_{\tau-}, \mathrm{d}y) \, \mathrm{d}\tau$$

$$\leq \frac{(1 - \mathrm{e}^{-2\sigma t}) V^2 \|f\|_{\mathrm{Lip}_d}^2}{2\sigma}.$$

By Kallenberg (1997), Lemma 23.19, the process given for any $\tau > 0$ by

$$(\exp\{\tau Z_s^f - b^{-2}\|f\|_{\mathrm{Lip}_d}^{-2}(\mathrm{e}^{\tau b\|f\|_{\mathrm{Lip}_d}} - \tau b\|f\|_{\mathrm{Lip}_d} - 1)\langle Z^f, Z^f\rangle_s\})_{0 \leq s \leq t}$$

is a $\mathbb{P}_x$-supermartingale with respect to $(\mathscr{F}_s)_{0 \leq s \leq t}$. Thus, using the two previous estimates, we get for any $\tau > 0$:

$$\mathbb{E}_x[\mathrm{e}^{\tau(f(X_t) - \mathbb{E}_x[f(X_t)])}] = \mathbb{E}_x[\mathrm{e}^{\tau Z_t^f}]$$

$$\leq \exp\left\{\frac{(1 - \mathrm{e}^{-2\sigma t})V^2}{2\sigma b^2}(\mathrm{e}^{\tau b\|f\|_{\mathrm{Lip}_d}} - \tau b\|f\|_{\mathrm{Lip}_d} - 1)\right\}.$$

To remove the boundedness assumption on $f$, consider the sequence of bounded functions $f_n := \max\{-n, \min\{f, n\}\}$ converging pointwise to $f$. Then it is routine to show that $(f_n)_{n \in \mathbb{N}}$ is uniformly integrable with respect to the probability measure $P_t(x, \cdot)$, which implies the $L^1$-convergence. Finally, since the functions $f_n$ are Lipschitz with a constant of at most $\|f\|_{\mathrm{Lip}_d}$ and $h$ is non-decreasing in its last variable, the use of Fatou's lemma achieves the proof. $\square$

Our present purpose is to extend to the multidimensional case the Laplace transform estimate (3.1) by using the method of tensorization.

Given $n \in \mathbb{N} \setminus \{0\}$, define $\mathrm{Lip}_{d_n}(\mathcal{X}^n)$ as the space of real Lipschitz functions $f$ on the product space $\mathcal{X}^n$, endowed with the seminorm

$$\|f\|_{\mathrm{Lip}_{d_n}} := \sup_{x \neq y} \frac{|f(x) - f(y)|}{d_n(x,y)} < +\infty,$$

where $d_n$ is the $\ell^1$-distance on $\mathcal{X}^n$ with respect to the metric $d$, that is, $d_n(x,y) := \sum_{i=1}^n d(x_i, y_i)$, $x, y \in \mathcal{X}^n$.

**Lemma 3.2.** *We assume that the hypothesis of Theorem 2.6 is fulfilled. Define the sample $X^n$ of the process $(X_t)_{t \geq 0}$ by $X^n = (X_{t_1}, \ldots, X_{t_n})$, $0 =: t_0 < t_1 < \cdots < t_n$ and let $f \in \mathrm{Lip}_{d_n}(\mathcal{X}^n)$. Then for any initial state $x \in \mathcal{X}$ and any $\tau > 0$, we have the multidimensional Laplace transform estimate:*

$$\mathbb{E}_x[e^{\tau(f(X^n) - \mathbb{E}_x[f(X^n)])}] \leq \exp\left\{\sum_{k=1}^n h(\tau, t_k - t_{k-1}, s_k b \|f\|_{\mathrm{Lip}_{d_n}})\right\}, \qquad (3.3)$$

*where the function $h$ is defined in Lemma 3.1 and $s_k := \sum_{l=k}^n e^{-\sigma(t_l - t_k)}$.*

**Proof.** Let $f_n := f$ and define for any $k = 1, \ldots, n-1$, the function $f_k$ on $\mathcal{X}^k$ by

$$f_k(x_1, \ldots, x_k) := \int_{\mathcal{X}^{n-k}} f(x_1, \ldots, x_n) P_{t_n - t_{n-1}}(x_{n-1}, dx_n) \cdots P_{t_{k+1} - t_k}(x_k, dx_{k+1})$$

$$= \int_{\mathcal{X}} f_{k+1}(x_1, \ldots, x_k, x_{k+1}) P_{t_{k+1} - t_k}(x_k, dx_{k+1}).$$

We divide the proof of Lemma 3.2 into two parts.

• Step 1: By a downward recursive argument on $k$, let us show first that the univariate function $x_k \mapsto f_k(*, x_k)$ is Lipschitz with respect to the metric $d$, with furthermore the inequality

$$\sup_{x_1, \ldots, x_{k-1} \in \mathcal{X}} \|f_k(x_1, \ldots, x_{k-1}, \cdot)\|_{\mathrm{Lip}_d} \leq s_k \|f\|_{\mathrm{Lip}_{d_n}}. \qquad (3.4)$$

Since $s_n = 1$, the property (3.4) is trivially true for $k = n$.

Assume now that (3.4) is satisfied for some $k \in \{2, \ldots, n\}$. First, letting $x_1, \ldots, x_{k-2}, y, z, x_k \in \mathcal{X}$, we have:

$$|f_k(x_1, \ldots, x_{k-2}, y, x_k) - f_k(x_1, \ldots, x_{k-2}, z, x_k)|$$

$$= \left| \int_{\mathcal{X}^{n-k}} f(x_1, \ldots, x_{k-2}, y, x_k, x_{k+1}, \ldots, x_n) P_{t_n - t_{n-1}}(x_{n-1}, dx_n) \cdots P_{t_{k+1} - t_k}(x_k, dx_{k+1}) \right.$$

$$\left. - \int_{\mathcal{X}^{n-k}} f(x_1, \ldots, x_{k-2}, z, x_k, x_{k+1}, \ldots, x_n) P_{t_n - t_{n-1}}(x_{n-1}, dx_n) \cdots P_{t_{k+1} - t_k}(x_k, dx_{k+1}) \right|$$

$$\leq \|f\|_{\mathrm{Lip}_{d_n}} d(y, z) \int_{\mathcal{X}^{n-k}} P_{t_n - t_{n-1}}(x_{n-1}, dx_n) \cdots P_{t_{k+1} - t_k}(x_k, dx_{k+1})$$

$$= \|f\|_{\mathrm{Lip}_{d_n}} d(y, z),$$

from which follows the inequality

$$\sup_{x_1,\ldots,x_{k-2},x_k \in \mathcal{X}} \|f_k(x_1,\ldots,x_{k-2},\cdot,x_k)\|_{\mathrm{Lip}_d} \leq \|f\|_{\mathrm{Lip}_{d_n}}. \tag{3.5}$$

Now, let us show that the property (3.4) is satisfied at the step $k-1$ with the help of (3.5). Let $x_1,\ldots,x_{k-2},y,z \in \mathcal{X}$. Using the contraction property (2.2) in the second inequality below,

$$|f_{k-1}(x_1,\ldots,x_{k-2},y) - f_{k-1}(x_1,\ldots,x_{k-2},z)|$$
$$\leq \left| \int_{\mathcal{X}} f_k(x_1,\ldots,x_{k-2},y,x_k)(P_{t_k-t_{k-1}}(y,\mathrm{d}x_k) - P_{t_k-t_{k-1}}(z,\mathrm{d}x_k)) \right|$$
$$+ \int_{\mathcal{X}} |f_k(x_1,\ldots,x_{k-2},y,x_k) - f_k(x_1,\ldots,x_{k-2},z,x_k)| P_{t_k-t_{k-1}}(z,\mathrm{d}x_k)$$
$$\leq \mathrm{e}^{-\sigma(t_k-t_{k-1})} \|f_k(x_1,\ldots,x_{k-2},y,\cdot)\|_{\mathrm{Lip}_d} d(y,z)$$
$$+ \int_{\mathcal{X}} \|f_k(x_1,\ldots,x_{k-2},\cdot,x_k)\|_{\mathrm{Lip}_d} d(y,z) P_{t_k-t_{k-1}}(z,\mathrm{d}x_k)$$
$$\leq (s_k \mathrm{e}^{-\sigma(t_k-t_{k-1})} + 1) \|f\|_{\mathrm{Lip}_{d_n}} d(y,z)$$
$$= s_{k-1} \|f\|_{\mathrm{Lip}_{d_n}} d(y,z),$$

where in the last inequality we used assumption (3.4) at the step $k$ together with (3.5). Therefore, we obtain the inequality

$$\|f_{k-1}(x_1,\ldots,x_{k-2},\cdot)\|_{\mathrm{Lip}_d} \leq s_{k-1} \|f\|_{\mathrm{Lip}_{d_n}},$$

and the parameters $x_1,\ldots,x_{k-2}$ being arbitrary, the property (3.4) is established at the step $k-1$, hence in full generality.

• Step 2: Proof of the Laplace transform estimate (3.3).

As before, let us show by a downward recursive argument on $k \in \{2,\ldots,n\}$ the following inequality:

$$\mathbb{E}_x[\mathrm{e}^{\tau f(X^n)}] \leq \exp\left\{ \sum_{i=k}^{n} h(\tau, t_i - t_{i-1}, bs_i \|f\|_{\mathrm{Lip}_{d_n}}) \right\}$$
$$\times \int_{\mathcal{X}^{k-1}} \mathrm{e}^{\tau f_{k-1}(x_1,\ldots,x_{k-1})} P_{t_{k-1}-t_{k-2}}(x_{k-2},\mathrm{d}x_{k-1}) \cdots P_{t_1}(x,\mathrm{d}x_1). \tag{3.6}$$

First let $k = n$. By the Markov property, we have

$$\mathbb{E}_x[\mathrm{e}^{\tau f(X^n)}]$$
$$= \int_{\mathcal{X}^n} \mathrm{e}^{\tau f_n(x_1,\ldots,x_n)} P_{t_n-t_{n-1}}(x_{n-1},\mathrm{d}x_n) \cdots P_{t_1}(x,\mathrm{d}x_1)$$

$$\leq \exp\{h(\tau, t_n - t_{n-1}, b\|f\|_{\mathrm{Lip}_{d_n}})\}$$
$$\times \int_{\mathcal{X}^{n-1}} e^{\tau f_{n-1}(x_1,\ldots,x_{n-1})} P_{t_{n-1}-t_{n-2}}(x_{n-2}, \mathrm{d}x_{n-1}) \cdots P_{t_1}(x, \mathrm{d}x_1),$$

where we used Lemma 3.1 with the univariate Lipschitz function $x_n \mapsto f_n(*, x_n)$ together with the inequality (3.4) since the function $h$ is non-decreasing in its last variable. Hence (3.6) is established in the case $k = n$.

Now assume that (3.6) is satisfied for some $k \in \{2, \ldots, n\}$. Using the same reasoning as above with the Lipschitz function $x_{k-1} \mapsto f_{k-1}(*, x_{k-1})$, we obtain

$$\mathbb{E}_x[e^{\tau f(X^n)}] \leq \exp\left\{\sum_{i=k}^n h(\tau, t_i - t_{i-1}, bs_i \|f\|_{\mathrm{Lip}_{d_n}})\right\}$$
$$\times \int_{\mathcal{X}^{k-1}} e^{\tau f_{k-1}(x_1,\ldots,x_{k-1})} P_{t_{k-1}-t_{k-2}}(x_{k-2}, \mathrm{d}x_{k-1}) \cdots P_{t_1}(x, \mathrm{d}x_1)$$
$$\leq \exp\left\{\sum_{i=k-1}^n h(\tau, t_i - t_{i-1}, bs_i \|f\|_{\mathrm{Lip}_{d_n}})\right\}$$
$$\times \int_{\mathcal{X}^{k-2}} e^{\tau f_{k-2}(x_1,\ldots,x_{k-2})} P_{t_{k-2}-t_{k-3}}(x_{k-3}, \mathrm{d}x_{k-2}) \cdots P_{t_1}(x, \mathrm{d}x_1)$$

so that the inequality (3.6) is satisfied at step $k-1$, hence in full generality. Finally, we obtain from (3.6) with $k=2$ the inequality

$$\mathbb{E}_x[e^{\tau f(X^n)}] \leq \exp\left\{\sum_{i=2}^n h(\tau, t_i - t_{i-1}, bs_i \|f\|_{\mathrm{Lip}_{d_n}})\right\} \int_{\mathcal{X}} e^{\tau f_1(x_1)} P_{t_1}(x, \mathrm{d}x_1)$$

and, using once again the same reasoning as before for the Lipschitz function $f_1$ entails the desired estimate (3.3). The proof of Lemma 3.2 is complete. $\square$

Now we are able to prove Theorem 2.6, with the help of Lemma 3.2.

**Proof of Theorem 2.6.** Define the sample $X^n = (X_{t_1}, \ldots, X_{t_n})$, where the sequence $t_k = kt/n$, $k = 0, \ldots, n$, is a regular subdivision of the time interval $[0, t]$. Since $\phi \in \mathrm{Lip}_d(\mathcal{X})$, the function $f$ given by $f(z_1, \ldots, z_n) := n^{-1} \sum_{k=1}^n \phi(z_k)$, $(z_1, \ldots, z_n) \in \mathcal{X}^n$, is Lipschitz on the product space $\mathcal{X}^n$ with respect to the $\ell^1$-metric $d_n$ and its Lipschitz seminorm satisfies $\|f\|_{\mathrm{Lip}_{d_n}} \leq n^{-1}\|\phi\|_{\mathrm{Lip}_d}$. Note that the function $h$ defined by (3.2) is non-decreasing in its last variable. Hence, since we have

$$\sup_{k=1,\ldots,n} s_k = \frac{1 - e^{-\sigma t}}{1 - e^{-\sigma t/n}},$$

the multidimensional Laplace transform estimate (3.3) of Lemma 3.2 implies the following upper bound:

$$\mathbb{E}_x[e^{\tau(f(X^n)-\mathbb{E}_x[f(X^n)])}] \leq \exp\left\{nh\left(\tau, \frac{t}{n}, \frac{b(1-e^{-\sigma t})\|\phi\|_{\mathrm{Lip}_d}}{n(1-e^{-\sigma t/n})}\right)\right\}, \qquad \tau > 0.$$

Therefore, by Chebyshev's inequality, we get for any $y > 0$:

$$\mathbb{P}_x(f(X^n) - \mathbb{E}_x[f(X^n)] \geq y)$$
$$\leq \inf_{\tau>0} e^{-\tau y} \mathbb{E}_x[e^{\tau(f(X^n)-\mathbb{E}_x[f(X^n)])}]$$
$$\leq e^{-(nV^2/(2b^2\sigma))(1-e^{-2\sigma t/n})g(2by\sigma(1-e^{-\sigma t/n})/(V^2(1-e^{-2\sigma t/n})(1-e^{-\sigma t})\|\phi\|_{\mathrm{Lip}_d}))}.$$

Applying also the same reasoning to the function $-f$ yields

$$\mathbb{P}_x(|f(X^n) - \mathbb{E}_x[f(X^n)]| \geq y)$$
$$\leq 2e^{-(nV^2/(2b^2\sigma))(1-e^{-2\sigma t/n})g(2by\sigma(1-e^{-\sigma t/n})/(V^2(1-e^{-2\sigma t/n})(1-e^{-\sigma t})\|\phi\|_{\mathrm{Lip}_d}))}. \tag{3.7}$$

Now, using the invariance property of the stationary distribution $\pi$ and the contraction property (2.2),

$$|\mathbb{E}_x[f(X^n)] - \pi(\phi)| = \left|\frac{1}{n}\sum_{k=1}^n \int_{\mathcal{X}} (P_{kt/n}\phi(x) - P_{kt/n}\phi(y))\pi(\mathrm{d}y)\right|$$
$$\leq \frac{1}{n}\sum_{k=1}^n e^{-\sigma kt/n} \|\phi\|_{\mathrm{Lip}_d} \int_{\mathcal{X}} d(x,y)\pi(\mathrm{d}y)$$
$$\leq \frac{(1-e^{-\sigma t})\|\phi\|_{\mathrm{Lip}_d}}{t\sigma} \int_{\mathcal{X}} d(x,y)\pi(\mathrm{d}y)$$
$$= M_t^x.$$

Hence the inequality (3.7) entails for any $y > 0$,

$$\mathbb{P}_x\left(\left|\frac{1}{n}\sum_{k=1}^n \phi(X_{kt/n}) - \pi(\phi)\right| \geq y + M_t^x\right) \leq 2e^{-A_n}, \tag{3.8}$$

where

$$A_n := \frac{nV^2}{2b^2\sigma}(1-e^{-2\sigma t/n})g\left(\frac{2by\sigma(1-e^{-\sigma t/n})}{V^2(1-e^{-2\sigma t/n})(1-e^{-\sigma t})\|\phi\|_{\mathrm{Lip}_d}}\right).$$

To finish the proof, note that since the process $(X_t)_{t\geq 0}$ is cadlag and the function $\phi$ is Lipschitz, the process $(\phi(X_t))_{t\geq 0}$ itself is cadlag so that the Riemann sum $n^{-1} \sum_{k=1}^{n} \phi(X_{kt/n})$ converges $\mathbb{P}_x$-a.s. to the empirical mean $t^{-1} \int_0^t \phi(X_s) \, \mathrm{d}s$. Therefore, using Fatou's lemma and the estimate (3.8), we obtain

$$\mathbb{P}_x\left(\left|\frac{1}{t}\int_0^t \phi(X_s)\,\mathrm{d}s - \pi(\phi)\right| \geq y + M_t^x\right) \leq \liminf_{n\to+\infty} \mathbb{P}_x\left(\left|\frac{1}{n}\sum_{k=1}^n \phi(X_{kt/n}) - \pi(\phi)\right| \geq y + M_t^x\right)$$

$$\leq \liminf_{n\to+\infty} 2\mathrm{e}^{-A_n}$$

$$= 2\mathrm{e}^{-(V^2 t/b^2) g((by\sigma)/(V^2(1-\mathrm{e}^{-\sigma t})\|\phi\|_{\mathrm{Lip}_d}))}.$$

The proof of Theorem 2.6 is established. □

## 4. Application to birth–death processes

The purpose of this final part is to apply Theorem 2.6 to birth–death processes. To do so, we compute the associated Wasserstein curvature with respect to a large class of metrics on $\mathbb{N}$. In particular, choosing suitably the metric with respect to the transition rates of the generator allows us to consider processes with non-necessarily bounded generators such as the classical $M/M/\infty$ queueing process.

Let $(X_t)_{t\geq 0}$ be a birth–death process on the state space $\mathcal{X} = \mathbb{N}$. This is a Markov process with a generator given for any function $f:\mathbb{N} \to \mathbb{R}$ by

$$\mathcal{L}f(x) = \lambda_x(f(x+1) - f(x)) + \nu_x(f(x-1) - f(x)), \qquad x \in \mathbb{N},$$

where the transition rates $\lambda$ and $\nu$ are positive with $\nu_0 = 0$, conditions ensuring the irreducibility of the process. Letting

$$\mu(0) := 1, \qquad \mu(x) := \frac{\lambda_0 \lambda_1 \cdots \lambda_{x-1}}{\nu_1 \nu_2 \cdots \nu_x}, \qquad x \geq 1,$$

we assume in the sequel that the process is ergodic, that is, it satisfies the following properties:

$$\sum_{x\geq 0} \mu(x) \sum_{y\geq x} \frac{1}{\mu(y)\lambda_y} = +\infty, \qquad C := \sum_{x\geq 0} \mu(x) < +\infty.$$

Then the stationary distribution of the process is $\pi(x) = \mu(x)/C$, $x \in \mathbb{N}$.

A fundamental example is the $M/M/\infty$ queue, also known as the birth–death process with immigration, which is an ergodic birth–death process $(X_t)_{t\geq 0}$ with an unbounded generator given by

$$\mathcal{L}f(x) = \lambda(f(x+1) - f(x)) + \nu x(f(x-1) - f(x)), \qquad x \in \mathbb{N},$$

where the parameters $\lambda$ and $\nu$ are positive. The associated stationary distribution is the Poisson measure $\mathscr{P}_\xi$ on $\mathbb{N}$ with parameter $\xi := \lambda/\nu$. Denote by $\mathscr{B}_{n,p}$ the binomial distribution with parameters $n \in \mathbb{N}$ and $p \in (0,1)$. Using the Mehler-type convolution formula given by Chafaï (2006):

$$\mathcal{L}(X_t | X_0 = x) = \mathscr{B}_{x, e^{-\nu t}} * \mathscr{P}_{\xi(1 - e^{-\nu t})}, \qquad t > 0,$$

we obtain by Chebyshev's inequality the following estimate, available for any $y > 0$:

$$\begin{aligned}
\mathbb{P}_x(X_t - \mathbb{E}_x[X_t] \geq y) &\leq \inf_{\tau > 0} e^{-\tau y} \mathbb{E}_x[e^{\tau(X_t - \mathbb{E}_x[X_t])}] \\
&\leq \inf_{\tau > 0} e^{-\tau y + \mathbb{E}_x[X_t](e^\tau - \tau - 1)} \\
&= \exp\left\{ y - (\mathbb{E}_x[X_t] + y) \log\left(1 + \frac{y}{\mathbb{E}_x[X_t]}\right) \right\},
\end{aligned}$$

where in the second inequality we used $\log(1+u) \leq u$, $u > 0$. Note that the latter Poisson-type deviation inequality is convenient for large time since we recover as $t$ tends to infinity the classical tail estimate for a centered Poisson random variable $X$ with intensity $\xi$:

$$\mathbb{P}(X - \mathbb{E}[X] \geq y) \leq \exp\left\{ y - (\xi + y) \log\left(1 + \frac{y}{\xi}\right) \right\}.$$

On the one hand, the $M/M/\infty$ queueing process is a discrete approximation of the Ornstein–Uhlenbeck process, whose stationary distribution is Gaussian. On the other hand, Remark 2.10 states that under the Bakry–Émery curvature criterion, the empirical mean of a Langevin-type process, which generalizes the Ornstein–Uhlenbeck process, satisfies a Gaussian deviation inequality. Hence it is natural, by comparison with the diffusion framework, to investigate Poisson-type tail estimates for the empirical mean of positively curved birth–death processes, since they generalize similarly the $M/M/\infty$ queueing process. However, if we consider the classical metric on $\mathbb{N}$, we are not able to apply Theorem 2.6 to processes with unbounded generators because, in this case, $V$ is infinite. Since the Wasserstein curvature strongly depends on the metric, the idea to overcome this difficulty is to carry the analysis with a Wasserstein curvature related to another metric on $\mathbb{N}$ that we choose suitably.

**Definition 4.1.** *Given a positive function $u$ on $\mathbb{N}$, define the metric $\delta : \mathbb{N} \times \mathbb{N} \to [0, +\infty)$ as*

$$\delta(x,y) := \left| \sum_{k=0}^{x-1} u_k - \sum_{k=0}^{y-1} u_k \right|, \qquad u_{-1} := 1.$$

Let us compute the Wasserstein curvature associated to this metric. To do so, we use the notion of coupling operators initiated by Chen (1986).

**Definition 4.2.** *An operator $\tilde{\mathcal{L}}$ acting on the space of real-valued functions on $\mathbb{N}^2$ is said to be a coupling of the generator $\mathcal{L}$ if it satisfies the two following properties:*

(i) *Marginality:*

$$\begin{cases} \tilde{\mathcal{L}} f_1(x,y) = \mathcal{L} f_1(x), \\ \tilde{\mathcal{L}} f_2(x,y) = \mathcal{L} f_2(y); \end{cases}$$

(ii) *Normality:* $\tilde{\mathcal{L}} h(x,x) = \mathcal{L} g(x)$.

*Here the two real-valued functions $f_1$ and $f_2$ on $\mathbb{N}$ are regarded as bivariate functions on $\mathbb{N}^2$, and $g$ is the univariate function $g(x) = h(x,x)$.*

Denote by $I$ the identity operator $I(f) = f$. Following Chen (1986), we introduce the classical coupling $\tilde{\mathcal{L}}$ by

$$\tilde{\mathcal{L}} f(x,y) = (\mathcal{L} \otimes I + I \otimes \mathcal{L}) f(x,y) 1_{\{x \neq y\}} + \mathcal{L} f(\cdot,\cdot)(x) 1_{\{x=y\}}, \qquad x, y \in \mathbb{N}.$$

Using the metric $\delta$, we have

$$\tilde{\mathcal{L}} \delta(x,y) = \begin{cases} \lambda_y u_y - \nu_y u_{y-1} - \lambda_x u_x + \nu_x u_{x-1}, & \text{if } x \leq y, \\ -\lambda_y u_y + \nu_y u_{y-1} + \lambda_x u_x - \nu_x u_{x-1}, & \text{otherwise.} \end{cases}$$

**Theorem 4.3.** *The Wasserstein curvature $\sigma_\delta$ with respect to the metric $\delta$ of the birth–death process $(X_t)_{t \geq 0}$ is given by the formula*

$$\sigma_\delta = \inf_{x \in \mathbb{N}} \left\{ \nu_{x+1} + \lambda_x - \nu_x \frac{u_{x-1}}{u_x} - \lambda_{x+1} \frac{u_{x+1}}{u_x} \right\}. \tag{4.1}$$

**Proof.** Denote $\alpha := \inf_{x \in \mathbb{N}} \{ \nu_{x+1} + \lambda_x - \nu_x \frac{u_{x-1}}{u_x} - \lambda_{x+1} \frac{u_{x+1}}{u_x} \}$ and assume first that $\sigma_\delta$ and $\alpha$ are not $-\infty$.

Consider on $\mathbb{N}$ the increasing Lipschitz function $f(x) = \sum_{k=0}^{x-1} u_k$ with Lipschitz seminorm $\|f\|_{\text{Lip}_\delta} = 1$. We have for any integers $x \leq y$ and any $t > 0$:

$$\frac{P_t f(y) - f(y)}{t} - \frac{P_t f(x) - f(x)}{t} = \frac{P_t f(y) - P_t f(x) - \delta(x,y)}{t}$$

$$\leq \frac{e^{-\sigma_\delta t} - 1}{t} \delta(x,y),$$

so that we obtain at the limit $t \to 0$:

$$\lambda_y u_y - \nu_y u_{y-1} - \lambda_x u_x + \nu_x u_{x-1} = \mathcal{L} f(y) - \mathcal{L} f(x) \leq -\sigma_\delta \delta(x,y).$$

Therefore, taking $y = x + 1$ and dividing by $u_x$ entail the inequality $\sigma_\delta \leq \alpha$.

On the other hand, we aim at proving that the Wasserstein curvature $\sigma_\delta$ is bounded below by $\alpha$. To do so, we use the coupling argument derived from the proof of Chen (1996), Theorem 1.1. Note that $\alpha$ rewrites as

$$\alpha = \inf_{x \in \mathbb{N}} \frac{-\tilde{\mathcal{L}} \delta(x, x+1)}{\delta(x, x+1)},$$

where $\tilde{\mathcal{L}}$ is the classical coupling operator defined above, so that we have

$$\tilde{\mathcal{L}}\delta(x, x+1) \leq -\alpha\delta(x, x+1), \qquad x \in \mathbb{N}.$$

As the following identities hold for any $x, y \in \mathbb{N}$ such that $x < y$:

$$\begin{cases} \tilde{\mathcal{L}}\delta(x,y) = \sum_{k=x}^{y-1} \tilde{\mathcal{L}}\delta(k, k+1), \\ \delta(x,y) = \sum_{k=x}^{y-1} \delta(k, k+1), \end{cases}$$

we get from the latter inequality and using the symmetry between $x$ and $y$ the inequality

$$\tilde{\mathcal{L}}\delta(x,y) \leq -\alpha\delta(x,y), \qquad x, y \in \mathbb{N}, \tag{4.2}$$

which ensures the contraction property (2.3), and so the desired estimate $\sigma_\delta \geq \alpha$. The proof is achieved in the finite case.

Finally, if at least $\sigma_\delta$ or $\alpha$ is $-\infty$, we are able to adapt the previous argument to show that both are actually infinite. $\square$

***Remark 4.4.*** Van Doorn (1985, 1987) proved that the spectral gap $\lambda_1$, which equals the so-called decay parameter in his papers, is actually the supremum of the Wasserstein curvatures given in Theorem 4.3 over the possible metrics $\delta$ defined in Definition 4.1. Later, such a result has been rediscovered by Chen (1996) with the coupling method emphasized in the proof above.

Once the metric $\delta$ has been introduced in full generality, let us introduce an assumption relating the weight $u$ and the transition rates of the generator of the birth–death process $(X_t)_{t\geq 0}$. We denote in the sequel $a \wedge b := \min\{a, b\}$.

***Assumption A.*** *There exist two constants $K > 0$ and $C > 0$ such that*

$$\left(\inf_{x \geq 0} \lambda_x\right) \wedge \left(\inf_{x \geq 1} \nu_x\right) \geq K \quad and \quad u_x \leq C\left(\frac{1}{\sqrt{\nu_{x+1}}} \wedge \frac{1}{\sqrt{\lambda_x}}\right), \qquad x \in \mathbb{N}.$$

Under Assumption A, we have a control on the metric $\delta$ as follows. The proofs are straightforward.

***Lemma 4.5.*** *Under Assumption A, the two inequalities below hold:*

(1) $\quad \delta(x,y) \leq \dfrac{C}{\sqrt{K}}|x-y|, \qquad x, y \in \mathbb{N};$

(2) $\quad \sup_{x \in \mathbb{N}} \lambda_x \delta(x, x+1)^2 + \nu_x \delta(x, x-1)^2 \leq 2C^2.$

***Remark 4.6.*** If at least one of the transition rates of the generator is unbounded, then the function $u$ vanishes at infinity so that the two metrics in Lemma 4.5(1) are not bi-Lipschitz equivalent. In particular, the identity function $f(x) = x$ is not Lipschitz on $\mathbb{N}$ with respect to the metric $\delta$.

Now we are able to state the following tail estimate for the empirical mean of the birth–death process $(X_t)_{t \geq 0}$.

**Corollary 4.7.** *Assume that the Wasserstein curvature $\sigma_\delta$ given by (4.1) is positive and that Assumption A is satisfied. Letting $\phi \in \mathrm{Lip}_\delta(\mathbb{N})$, for any initial state $x \in \mathbb{N}$, any $t > 0$ and any $y > 0$, we have the following Poisson-type deviation inequality:*

$$\mathbb{P}_x\left(\left|\frac{1}{t}\int_0^t \phi(X_s)\,\mathrm{d}s - \pi(\phi)\right| \geq y + M_t^x\right) \leq 2\mathrm{e}^{-2Ktg((y\sigma_\delta)/(2\sqrt{K}C(1-\mathrm{e}^{-\sigma_\delta t})\|\phi\|_{\mathrm{Lip}_\delta}))}, \quad (4.3)$$

*where $M_t^x := \sigma_\delta^{-1} t^{-1}(1 - \mathrm{e}^{-\sigma_\delta t})\|\phi\|_{\mathrm{Lip}_\delta} \sum_{z \in \mathbb{N}} \delta(x,z)\pi(z)$ and $g$ is the function given in (2.5).*

**Proof.** Using Lemma 4.5, we get the result by applying Theorem 2.6 with $b = C/\sqrt{K}$ and $V^2 = 2C^2$. □

***Remark 4.8.*** The Poisson-type deviation inequality (4.3) is comparable to that obtained recently by Liu and Ma (2009) by using martingale techniques together with the so-called Lipschitz spectral gap. We mention, however, that there is a one-to-one correspondence between this object and the Wasserstein curvature according to the variational formulas given by Chen (1996), Theorem 1.1.

To finish this work, let us return to the case of the $M/M/\infty$ queueing process. For the sake of simplicity, we assume in the sequel that the intensity $\xi$ of the process equals 1. Choosing $u_x := (x+1)^{-1/2}$, $x \in \mathbb{N}$, in the definition of the metric $\delta$, a brief computation shows that the Wasserstein curvature $\sigma_\delta$ equals $\nu/2$, which is half of the exact curvature $\nu$ given by Chafaï (2006). Moreover, the transition rates of the generator satisfy Assumption A with $C = \sqrt{K} = \sqrt{\nu}$. Hence, Corollary 4.7 entails for any Lipschitz function $\phi \in \mathrm{Lip}_\delta(\mathbb{N})$, any $t > 0$, any initial state $x \in \mathbb{N}$ and any $y > 0$,

$$\mathbb{P}_x\left(\left|\frac{1}{t}\int_0^t \phi(X_s)\,\mathrm{d}s - \mathscr{P}_1(\phi)\right| \geq y + M_t^x\right) \leq 2\mathrm{e}^{-2\nu tg(y/(4(1-\mathrm{e}^{-\nu t/2})\|\phi\|_{\mathrm{Lip}_\delta}))}.$$

***Remark 4.9.*** An inequality such as the one above allows us to consider unbounded functions $\phi$ as, for instance, the square root function, which is Lipschitz with respect to the metric $\delta$. However, as noted in Remark 4.6, the price to pay is to require $\phi \in \mathrm{Lip}_\delta(\mathbb{N})$, which unfortunately excludes the identity function since the generator is unbounded. Hence we conjecture that in the case of the $M/M/\infty$ queueing process, the deviation of the empirical mean of Lipschitz functions with respect to the classical metric is of the Poisson type. See also the recent work of Guillin *et al.* (2009) for an approach to this problem through transportation-information inequalities.